\documentclass[fleqn,11pt,a4paper]{article}

\usepackage{a4wide}
\usepackage{graphicx}

\usepackage{color}
\usepackage[latin1]{inputenc}

\usepackage{amsfonts}
\usepackage{amsmath}
\usepackage{amssymb}
\usepackage{amsthm}
\usepackage{enumerate}
\usepackage{setspace}
\newtheorem{Theorem}{Theorem}

\newtheorem{Lemma}[Theorem]{Lemma}
\newtheorem{Corollary}[Theorem]{Corollary}
\newtheorem{Proposition}[Theorem]{Proposition}
 \theoremstyle{definition}
\newtheorem{Definition}[Theorem]{Definition}
\newtheorem{Remark}[Theorem]{Remark}
\newtheorem{Example}[Theorem]{Example}

\everymath{\displaystyle}

\newcommand{\+}[1]{{\cal {#1}} }

\def\Wk{{\+A}^k}

\def\vep{{\varepsilon}}

\newcommand{\R}{{\mathbb R}}

\begin{document}

\title{Perfect Necklaces}

\author{\begin{tabular}{cccc}
Nicol\'as Alvarez  &   Ver\'onica Becher &
Pablo A. Ferrari    &   Sergio A. Yuhjtman
\end{tabular}}

\maketitle

\begin{abstract}
We introduce a  variant of de~Bruijn words that we call perfect necklaces.
Fix a finite  alphabet. Recall that a word  is a finite sequence of symbols in the alphabet
and a  circular word, or necklace, is the equivalence class of a word under  rotations.
For  positive integers  $k$ and $n$, we call a necklace  $(k,n)$-perfect  if 
each word of length $k$ occurs exactly~$n$ times  at  positions which
are different modulo~$n$ for any convention on the starting point.
We call a necklace perfect if it is $(k,k)$-perfect for some~$k$.
We prove that  every arithmetic sequence with difference coprime with the alphabet size induces a perfect necklace.
In particular,  the concatenation of all words of the same length in lexicographic order  yields a perfect necklace.
For each $k$ and $n$,  we give a closed formula for the  number of $(k,n)$-perfect necklaces.
Finally, we prove that every infinite periodic sequence whose period coincides 
with some $(k,n)$-perfect necklace for any $n$, passes all statistical tests of size up to~$k$, 
but not all  larger tests. This~last theorem motivated this work.
 \end{abstract}
\bigskip


\bigskip


\section{Introduction}

Fix a finite alphabet $\+A$ and write  $|\+A|$ for its  cardinality.
A word  is a finite sequence of symbols in the alphabet.
A circular word, or necklace, is the equivalence class of a word under rotations.
In this note we introduce  {\em perfect necklaces}:

\begin{Definition}
A necklace   is \emph{$(k,n)$-perfect} if it has length~$n |\+A|^k$
and each word of length $k$ occurs exactly~$n$ times  at  positions which are different modulo~$n$ 
for any convention on the starting point.
A necklace is {\em perfect} if it is $(k,k)$-perfect for some~$k$.
\end{Definition}
%

Perfect necklaces are a variant  of the celebrated de~Bruijn necklaces~\cite{debruijn}.
Recall that a  de~Bruijn necklace  of order~$k$ in alphabet~$\+A$  has  length $|\+A|^k$ 
and each word of length~$k$ occurs in it exactly once. 
Thus, our $(k,1)$-perfect necklaces coincide with the de Bruijn necklaces of order~$k$.
For a supreme presentation of de~Bruijn necklaces, including a historic account of their discovery and rediscovery,
 see~\cite{berstel2007}. Observe that a necklace of length $k|\+A|^k$ admits $k$ possible decompositions  
into $|\+A|^k$ consecutive (non-overlapping) words of length $k$. Hence, a necklace is $(k,k)$-perfect if and only if 
it has length $k |A|^k$ and each word of length $k$ occurs exactly once in each of the $k$ possible decompositions.
 
For each $k$ and $n$,  we give a characterization of $(k,n)$-perfect necklaces 
in terms of  Eulerian circuits in appropriate graphs (Corollary \ref{cor:ldifferent}).
We give  a closed formula for the number of  $(k,n)$-perfect necklaces (Theorem~\ref{thm:count}).
These are the most elaborate results in this work.

We show that each arithmetic sequence with difference coprime with the alphabet size induces a perfect necklace 
(Theorem~\ref{thm:arithmetic}).
In particular, the concatenation of  all words of the same length in lexicographic order
yields a perfect necklace (Corollary~\ref{cor:ordered}). This provides a gracious instance of a perfect necklace
for any word length.

As far as we know, David Champernowne~\cite{Champernowne1933}
was the first to consider  combinatorial properties in the concatenation of  all words of the same length
in lexicographic order. He used them in his construction of a real number  normal to base~$10$, 
a property defined by  \'Emile Borel~\cite{Borel1909}.
He worked with alphabet $\+A=\{0, 1, \ldots, 9\}$ and for each~$k$, 
he  bounded the number of  occurrences of each word of length up to~$k$ in 
 the concatenation of  all words of length $k$ in lexicographic order.
But Champernowne missed  that each word of length~$k$ occurs in this sequence exactly~$k$ times, 
once in each of the~$k$ different~shifts.

\section{Perfect  necklaces}

\paragraph{Notation.} 
We write  ${\+A}^*$  for the set of all words, and $\+A^k$ for the set of all words of length~$k$.
The length of a word $w$ is denoted with $|w|$ and the  positions in $w$ are numbered from $0$ to $|w|-1$.
We write  $w(i)$ to denote the symbol in the $i$-th position of $w$.
Let $\theta: \+A^* \to \+A^*$ be the \emph{shift} operator, such that for each  position~$i$,
$(\theta w)(i)=w((i+1) \mod |s|)).$
That is, the shift operator  is defined with the   convention of periodicity.
With~$\theta^n$ we denote  the application of the shift $n$~times to the right, and 
with~$\theta^{-n}$, $n$ times to the left.
As already stated, a necklace is the equivalence class of a word under rotations.
To denote a necklace we write $[w]$ where $w$~is any of the words in the equivalence class. 
For example, if $\+A=\{0,1\}$,

$[000]$  contains  a single word $000$, because for every $n$, $\theta^n(000)=000$.

$[110]$ contains three words $\theta^0(110)=110$, $\theta^1(110)=101$ and $\theta^2(110)=011$.


\begin{Example}  Let $\+A=\{0,1\}$. We add spaces in the examples  just for readability.

\noindent
For words of length $2$ there are just two perfect necklaces:

$[00\ 01\ 10\ 11]$, 

$[00\ 10\ 01\ 11]$.

\noindent
This is a perfect necklace for word length $3$:

$[000\;110\;101\;111\;001\;010\;011\;100]$.

\noindent
The following are not perfect,

$[00\;01\;11\;10]$,

$[000\;101 \;110\;111\;010\;001\;011\;100]$.

\noindent
The so-called \emph{Gray numbers} are not perfect, for instance,  
$[000\;001\;011\;010\;110\;111\;101\;100].$
\end{Example}

\subsection{Each ordered necklace is perfect}

\begin{Definition}
For an ordered alphabet $\+A$  and  a positive integer $k$, the $k$-\emph{ordered necklace}  
has length~$k|\+A|^k $ 
and it  is obtained by the concatenation of all words  of length $k$ in lexicographic order.
\end{Definition}

\noindent
For $\+A=\{0,1\}$ the following are the ordered necklaces for $k$ equal to $1$, $2$ and~$3$ respectively:

$ [0 1],$  

$ [00\ 01\ 10\ 11]$,

$ [000\ 001\ 010\ 011\  100\ 101\  110\  111]$.

We will prove that for every word length,   the ordered necklace is perfect.
 We say that a bijection $\sigma: \Wk\to \Wk$ is a {\em cycle} if for each $w\in\Wk$
the set $\{\sigma^j(w): 0\leq j < |\+A|^k \}$ equals $\Wk$.
For a word $w$ we write 
$w(i\ldots j)$ to denote the subsequence of $w$ from position $i$ to $j$.

\begin{Lemma}\label{lema:reescritura}
Let $\+A$ be a finite alphabet, 
$\sigma: \Wk \to \Wk$  a cycle and $v$ any word in $\+A^k$.
Let  $s=\sigma^0(v) \sigma^1(v) \ldots \sigma^{|\+A|^k-1}(v)$.
The necklace  $[s]$ is perfect if and only if 
for every $\ell$ such that $0 \leq \ell < k$, for every  $x \in \+A^\ell$ and every $y \in \+A^{k-\ell}$,
 there is a unique  $w \in \Wk$ 
such that $w(k-\ell\ldots k-1)=x$ and $(\sigma(w))(0\ldots k-\ell-1)=y$.
 \end{Lemma}
\begin{proof} 
 Assume  $[s]$ is $(k,k)$-perfect. 
 Take $\ell$ such that $0 \leq \ell < k$, $x \in {\+A}^\ell$ and $y \in {\+A}^{k-\ell}$. 
 Consider  $\theta^{-\ell} s$, the $-\ell^{th}$ shift of~$s$.
  Since $[s]$ is $(k,k)$-perfect, $xy$ occurs exactly once in the decomposition of 
 $\theta^{-\ell} s$ in consecutive words of length~$k$.
 Thus, there is a  unique word $w$ in the  decomposition of $s$ in consecutive words of length~$k$  whose last~$\ell$ 
symbols are equal to~$x$ and whose first~$k-\ell$ symbols are equal to~$y$.
 Conversely, suppose $[s]$ is not $(k,k)$-perfect. 
Then, there is some $\ell$, $0\leq\ell<k$, such that the decomposition 
of  $\theta^{-\ell}(s)$ contains two  equal words of length~$k$.
This contradicts that for every $x\in\+A^\ell$ and every $y\in \+A^{k-\ell}$, there is a unique $w\in \+A^k$ such that 
$w(k-\ell\ldots k-1)=x$ and $(\sigma(w))(0\ldots k-\ell-1)=y$.
 \end{proof}

\begin{Theorem}\label{thm:arithmetic}
 Consider the alphabet $\+A = \{0,..,b-1\}$ where $b$ is an integer greater than or equal to $2$,
 a word length~$k$ and a positive integer $r$  coprime with $b$.
 Identify the elements of $\+A^k$ with  the set of integers modulo $b^k$
 according to representation in base $b$.
 Define the word of length $kb^k$ by the juxtaposition of the elements of $\+A^k$ corresponding 
 to the arithmetic sequence $0, r, 2r, \dots, (b^k-1)r$. Then the associated necklace is perfect.
\end{Theorem}

\begin{proof} 
Since $r$ is coprime with $b$, the addition of $r$ defines a cycle $\sigma:\+A^k \to \+A^k$. 
We must check that  it satisfies the condition in Lemma~\ref{lema:reescritura}. 
For any $w$ such that $w(k-\ell\ldots k-1)=x$ we have $\sigma(w)(k-\ell\ldots k-1)=\tilde x$, 
where abusing notation $\tilde x=x+r\mod b^\ell$. 
Since the word $y\tilde x$ appears only one time in the cycle, 
this fixes a unique $w=\sigma^{-1}(y\tilde x)$ with $w(k-\ell\ldots k-1)=x$ and $(\sigma(w))(0\ldots k-\ell-1)=y$.  
\end{proof}

\begin{Corollary}\label{cor:ordered} 
For an ordered alphabet $\+A$  and word length~$k$, the $k$-ordered necklace is perfect.
\end{Corollary}

\begin{proof}
Take $r=1$ in Theorem \ref{thm:arithmetic}.
\end{proof}

The following proposition is immediate, so we state it without proof.

\begin{Proposition}  \label{symmetries}\samepage
The following operators $\phi:\+A^*\to \+A^*$  are well defined on necklaces and  preserve perfection.
That is, for every~$k$ and $n$ and for every $s\in \+A^*$, 
if $[s]$  is $(k,n)$-perfect then   $[\phi s]$ is~$(k,n)$-perfect.
\begin{enumerate}
\item The digit permutation operator defined by 
$\phi(x_0 \dots x_{kb^k-1})=(\pi x_0 \dots \pi x_{kb^k-1})$ for any permutation $\pi:\+A\to\+A$.
\samepage
\item The reflection operator $\phi(x_0 \dots x_{kb^k-1})=(x_{kb^k-1} \dots x_0)$.
\end{enumerate}
\end{Proposition}

\section{Characterizing and counting perfect necklaces}

To characterize and count $(k,n)$-perfect necklaces in alphabet $\+A$ we consider 
Eulerian circuits in an appropriate directed graph, defined from $\+A$, $k$ and $n$. 
Recall that an Eulerian circuit in a graph is a path that uses all edges exactly once.
A thorough presentation of the material on graphs that we use in this section can 
be read in the monographs~\cite{harary,tutte,cvetovic}.
For the material on combinatorics on words see  the books~\cite{lothaire1997,lothaire2002}.

We write $m|n$  when $m$ divides $n$ and we write $\gcd(m,n)$ for the maximum common divisor between $m$ and~$n$.

\begin{Definition}
Let  $\+A$ be an  alphabet with cardinality $b$, 
let $s$ be  a word length 
and let $n$ be a positive integer. 
We define the {\em astute graph} $G_{s,n}$ 
as the  directed graph, 
with $n b^s$ nodes, 
each node is a pair $(u,v)$, where $u$ is in $\+A^s$ 
and  $v$ is a number between  $0$ and  $n-1$.
There is an edge from $(u,v)$ to  $(u',v')$ 
if the last  $s-1$ symbols from  $u$ coincide with the first  $s-1$ symbols from~$u'$
and $(v+1) \mod n = v'$.
Observe that  $G_{s,n}$ is  strongly regular  (all nodes have in-degree  and out-degree equal to $b$) and 
it is strongly connected (there is a path from every node to every other node).
\end{Definition}

\begin{Remark}
For any alphabet size, the astute graph $G_{k-1,1}$ coincides with a de~Bruijn graph of words of length $k-1$; 
hence, the Eulerian circuits in $G_{k-1,1}$ yield exactly the de~Bruijn~necklaces of order~$k$.
\end{Remark}

Although each Eulerian circuit in the astute graph $G_{k-1,n}$  gives one $(k,n)$-perfect necklace, 
each $(k,n)$-perfect necklace can come from several Eulerian circuits in this graph.

\subsection{From perfect necklaces to Eulerian circuits}

Hereafter, we assume an alphabet $\+A$ and  we write $b$ for its~cardinality. 

\begin{Definition}
For a necklace of length $\ell$, $[a_0, a_2, \ldots a_{\ell-1}]$,
we define its {\em period}  as the minimum integer $L$ such that for every 
non-negative integer $j$,
$a_{j \mod \ell}=a_{(j+L) \mod \ell}$. 
Notice that  the period $L$ always exists, and necessarily $L|\ell$. 
If the period coincides with the length we say the necklace is~{\em irreducible}.
\end{Definition}

\begin{Definition}\label{def:d}
Let $m, n$ be positive integers. We define $d_{m,n} = \prod p_i^{\alpha_i}$ where $\{ p_i \}$ is the set of primes
that divide $m$, and $\alpha_i$ is the exponent of $p_i$ in the factorization of~$n$.
\end{Definition}

\begin{Proposition}
The period $L$ of a $(k,n)$-perfect necklace satisfies the following: 
\begin{enumerate}
\samepage
\item  $L=jb^k$ for $j|n$.

\item $d_{b,n} | j$.

\item The corresponding irreducible necklace of length $L = jb^k$ is $(k,j)$-perfect.
\end{enumerate}
\end{Proposition}

\begin{proof} 
Let $ [s]$  be $(k,n)$-perfect, with  $s= a_0 \ldots a_{nb^k-1}$.

1.  Since $[s]$ has length $n b^k$, 
we know $L|nb^k$.
Let's verify that  $b^k|L$. 
Since $[s]$ has period~$L$,   $[a_0 \ldots a_{L-1}]$  is a necklace where 
all words of length~$k$ occur the same number of times.
Otherwise, it would be impossible that they occur the same number of times in~$[s]$.
If each word of length~$k$ occurs $j$ times in $[a_0 \ldots a_{L-1}]$, 
then~$L = jb^k$. 
Since $jb^k|nb^k$, we conclude~$ j|n$.

2.  The word $ a_0 \ldots a_{k-1}$  occurs at position $0$ in~$s$
but also at positions $L, 2L, \ldots, (n/j -1)L$.
These positions are of the form $ qjb^k $ where $0 \leq q < n/j$.
These numbers must have  pairwise different congruences modulo~$n$.
Equivalently, the $n/j$ numbers of the form  $rb^k$, where $0 \leq q < n/j$,
are all pairwise different  modulo~$n$. . 
This last condition holds exactly when  $\gcd(b^k, n/j) = 1$,
which in turn is equivalent to  $\gcd(b, n/j) = 1$, 
which is equivalent to $d_{b,n} | j$.

3.  As argued in Point 1, in the necklace $[a_0 \ldots a_{L-1}]$ 
every word of length $k$ occurs the same number of times. 
If the positions of two occurrences of a given word were equal modulo~$j$ then they 
would be equal modulo~$n$, 
but this is impossible because $[s]$ is $(k,n)$-perfect.
\end{proof}

\begin{Proposition} 
Let $N$ be a $(k,j)$-perfect necklace.
If $n$ is such that  $d_{b,n} | j | n$ then the  necklace of length $n b^k$  
obtained by repeating $N$ exactly $n/j$ times is   $(k,n)$-perfect. 
\end{Proposition}

\begin{proof}
  Let $\tilde N$ be obtained by repeating $N$ exactly $n/j$ times.
  Then each word of length $k$ occurs in $\tilde N$ exactly $j \times  n / j = n $ times. 
  Take a word $w$ of length $k$ and let   $q_1, \ldots, q_j$, each between $0$ and $j b^k -1$,
  be the positions of the occurrences of $w$ in $N$ for some convention on the starting point.
  Then, $w$ occurs in  $\tilde N$  at positions  $q_i + jb^k t$, where $0 \leq t < n/j$. 
  Assume $q_{i_1} + jb^kt_1 \equiv q_{i_2} + jb^kt_2 \ (\mbox{mod } n)$. 
 Taking modulo~$j$ we conclude $i_1 = i_2$ because  $N$ is   $(k,j)$-perfect.  
  Then we have 
  $b^k t_1 \equiv b^k t_2 \ (\mbox{mod} \ \ n/j)$.
  Since $d_{b,n} | j$ we have $\gcd(b, n/j)=1$, so 
  $t_1 \equiv t_2 \ (\mbox{mod } n/j)$, which implies $t_1 = t_2$.
\end{proof}

\begin{Corollary} \label{cor:ldifferent}
  Assume an alphabet of $b$ symbols, with $b\geq 2$. Let $k$ and $n$ be  positive integers.  
  An Eulerian circuit in the astute graph $G_{k-1,n}$ induces a $(k,n)$-perfect necklace. 
  Each  $(k,n)$-perfect  necklace of period $jb^k$ corresponds  to  $j$ different eulerian circuits in $G_{k-1,j}$.
  Therefore, the number      of Eulerian circuits in  the astute graph~$G_{k-1,n}$  is
  \[
 e(n) = \sum_{d_{b,n} | j | n} j \ p(j),
\]
  where 
$p(j)$ is the number of irreducible   $(k,j)$-perfect necklaces. 
 \end{Corollary}

\subsection{The number of Eulerian circuits in the astute graphs}

Let $G$ be a directed graph with $n$ nodes. 
The adjacency matrix of a graph $G$  is the matrix
 $A(G)=(a_{i,j})_{i,j=1}^{n}$ where $a_{i,j}$ is the number of edges between node~$i$ and node~$j$. 
The characteristic polynomial~\cite{cvetovic} of a graph $G$ is defined as 
\[
{\+P}(G;x) = \text{determinant} (x I- A(G)),
\]
where $I$ is the identity matrix of dimension  $n\times n$.

The BEST  theorem (for the authors  Bruijn, van Aardenne-Ehrenfest, Smith and Tutte)  gives 
a product formula for the number of Eulerian circuits in directed  graphs.

\begin{Lemma}[BEST Theorem~\cite{harary}]\label{lemma:BEST}
 Let $G$ be regular  connected graph  with $n$ nodes.
 Let $v$ be a node of $G$ and let $r(G)$ be the number of spanning trees oriented towards~$v$. 
The number of Eulerian circuits in $G$ is  
\[
r(G) \cdot \prod_{v=1}^n (degree(v)-1)!
\]
\end{Lemma}

\begin{Lemma}[Hutschenreurther, Proposition 1.4~\cite{cvetovic}]\label{lemma:c}
Let $G$ be  a regular multigraph  with  $n$ nodes and  degree $b$.
For any of its nodes, the number of spanning trees  $r(G)$ oriented to it  is 
\[
r(G)=\frac{1}{n}\frac{\partial}{\partial x}{\+P}(G;x) |_{x=b}. 
\]
where $\frac{\partial}{\partial x}$ is the derivative with respect to $x$.
\end{Lemma}

Given a graph $G$, its line-graph  $\Gamma(G)$ is 
a graph such that
each node of $\Gamma(G)$ represents an edge of $G$; and
two nodes of $\Gamma(G)$ are adjacent if and only if their corresponding edges share a common node in  in $G$.
\begin{Lemma}[\cite{cvetovic}]
For any directed  graph $G$, regular and connected, 
\[
{\+P}(\Gamma(G); x)= x^{m-n} {\+P}(G;x),
\]
where $\Gamma(G)$ is the line-graph of $G$, $m$ is the number of edges of~$G$ and~$n$ is the number of nodes of~$G$.
\end{Lemma}

In the next lemma we write $\lambda$ for the empty word, namely the unique word in $\+A^0$.

\begin{Lemma}
Let $b$ be any alphabet size, $k$~be a word length, 
and  $j$ be an integer such that $gcd(b,k)|j|k$.
Let $G_{0,j}$ be the  graph with the set of nodes $\{(\lambda,0), (\lambda,1), \ldots(\lambda,j-1)\}$, 
with  $b$~edges from $(\lambda,i)$ to $(\lambda,  i+1\mod j)$.
Then,
${\+P}(G_{0,j};x)=   x^j - b^j$.
\end{Lemma}
\begin{proof}
\mbox{It is easy to check that ${\+P}(G_{0,j};x)=det(x I- A(G_{0,j} ))$, which is equal to~$x^j-b^j$.}
\end{proof}

\begin{Lemma} \label{lemma:eulerian}
Assume  an alphabet of  $b$ symbols with $b\geq 2$.
Let $k$ be a  word length and $j$ be a positive integer such that $\gcd(b,k)|j|k$.
The number  of Eulerian circuits in the astute graph $G_{k-1,j}$ is
$ (b!)^{j b^{k-1}}b^{- k}.$
\end{Lemma}

\begin{proof}
We write  $\Gamma(G)$ to denote the line graph of~$G$.
Notice that for every positive $s$ and for every~$j$, 
$G_{s,j} = \Gamma(G_{s-1,j}).$ In this proof the value $j$ will remain fixed.

Since ${G_{k-1,j}}$  has $j b^{k-1}$ nodes, each with  in-degree~$b$ (also out-degree~$b$), 
by Lemma~\ref{lemma:BEST} the number of Eulerian circuits in $G_{k-1,j}$ is
\[
r(G_{k-1,j}) \cdot \prod_{v=1}^{j b^{k-1}} (degree(v)-1)! = r(G_{k-1,j})\dot  (b-1)!^{j b^{k-1}}.
\]
The rest of the proof  is to determine $r(G_{k-1,j})$ using Lemma~\ref{lemma:c}.
\begin{align*}
{\+P}(G_{k-1,j};x )
&=
 {\+P}(\Gamma(G_{k-2,j}) ;x)
\\
&=  x^{b^{k-1}j -b^{k-2}j} {\+ P}(G_{k-2,j};x)
\\
&= x^{j( b^{k-1}-b^{k-2})} {\+P}(\Gamma(G_{k-3,j});x)
\\
&= x^{j(b^{k-1} -b^{k-2})} x^{j(b^{k-2}- b^{k-3})} {\+P}(G_{k-3,j};x) 
\\
&=x^{j (b^{k-1} - b^{k-3})} {\+P}(G_{k-3,j};x)
\\
& = \ldots 
\\
&=  x^{j (b^{k-1} - b^{0})}{\+P}(G_{0,j};x)
\\
&=  x^{j (b^{k-1} - 1)} (x^j - b^j).
\\
\frac{\partial}{\partial x}{\+P}(G_{k-1,j} ;x) 
&=  \frac{\partial}{\partial x}   x^{j (b^{k-1} - 1)} (x^j - b^j )
\\
&=     (j b^{k-1} - j) x^{j b^{k-1} - j -1}  (x^j - b^j ) +  x^{j b^{k-1} - j}  j x^{j-1}.
\\
\frac{\partial}{\partial x}{\+P}(G_{k-1,j} ;x)  |_{ x=b }\ \
 &=    b^{j b^{k-1} -j}  j  b^{j-1}.
\end{align*}
Finally, by Lemma~\ref{lemma:c},
\begin{align*}
r(G_{k-1,j})
=& \frac{1}{j b^{k-1}} \frac{\partial}{\partial x} {\+P}(G_{k-1,j} ;x)|_{ x=b} 
=   \frac{1}{j b^{k-1}}  b^{j b^{k-1} -j}  j  b^{j-1} = b^{j b^{k-1}  -k }.
\end{align*}
Hence, the total number Eulerian circuits in $G_{k-1,j}$ is
\[
b^{j b^{k-1} -k} ((b-1)!)^{j b^{k-1}} = b!^{j b^{k-1}}  b^{-k}.\qedhere
\]
\end{proof}

\subsection{The number of perfect necklaces}

Recall that by Definition \ref{def:d}, 
$d_{b,n} = \prod p_i^{\alpha_i}$,  where  $\{ p_i \}$  is the set of primes that divide both $b$~and~$n$, 
and $\alpha_i$ is  the exponent of $p_i$ in the factorization of~$n$. 
The Euler totient function $\varphi(n)$ counts the positive integers less than or equal to $n$ that are relatively prime to $n$.

\begin{Theorem}\label{thm:count}
Assume  an alphabet of~$b$ symbols, with $b\geq 2$. Let $k$ and $n$ be  positive  integers.
The~number  of $(k,n)$-perfect necklaces is 
 \[
 \frac{1}{n} \sum_{d_{b,n}|j|n} e(j) \varphi(n/j) 
\]
\samepage
where 
$e(j)= (b!)^{j b^{k-1}}b^{- k}$ is the number  of Eulerian circuits in  graph $G_{k-1,j}$
and $\varphi$ is Euler's totient function.
\end{Theorem}

\begin{proof}
Let $p(j)$ be the number of irreducible $(k,j)$-perfect necklaces. 
Then, the number of $(k,n)$-perfect necklaces is
\[
\sum_{d_{b,n}|j|n} p(j).
\]
Let $e(j)$ be the number of Eulerian circuits in the astute graph $G_{k-1,j}$.
By  Corollary~\ref{cor:ldifferent}, for each  $j$ such that $d_{b,n}|j|n$,
\[
e(j) = \sum_{d_{b,n}|\ell|j} \ell\  p(\ell).
\]
Notice that $d_{b,n} = d_{b,j}$.
For a lighter notation, in the rest of the proof we abbreviate $d_{b,n}$ as just $d$.
Then, writing each such $j$ as  a multiple of~$d$, we obtain that for each~$m$ such that~$md|n$,
\[
e(m d) = \sum_{i|m} i d \ p(i d).
\]
Let $ g(m) = e(m d)$  and  $ f(m) =  p(m d)\ m d$. 
Wrting $\mu$ for  the  M\"obius function   we obtain
\begin{align*}
f(m) =& \sum_{i|m }\mu(m/i) \ g(i).
\\
p(md) \ m d=&\sum_{i|m} \mu(m/i) \ e(i d).
\\
p(md) =&\frac{1}{md}\sum_{i|m} \mu(m/i)\ e(i d).
\\
\sum_{d|j|n} p(j) =&  \sum_{m|n/d} \frac{1}{md} \sum_{i|m} \mu(m/i)\ e(i d)
\\
  =& \sum_{i|n/d}  e(i d) \sum_{i|m|n/d} \frac{1}{md} \ \mu(m/i)\
\\
 =& \sum_{d|j|n} e(j)\sum_{j|q|n} \frac{1}{q} \  \mu(q/j).
\end{align*}
Applying the M\"obius inversion,
\[
\sum_{j|q|n} \frac{1}{q} \ \mu(q/j) = 
\sum_{r|n/j} \frac{1}{jr} \ \mu(r) = 
\frac{1}{n} \sum_{r|n/j} \frac{n/j}{r} \ \mu(r)
= \frac{1}{n} \ \varphi(n/j).
\]
We have used the identity $\varphi(m) = \sum_{r|m} \frac{m}{r} \mu(r)$, which is simply the inversion of $m = \sum_{r|m} \varphi(r)$.
 By Lemma~\ref{lemma:eulerian},
the number $e(j)$ of Eulerian circuits in the astute graph $G_{k-1,j}$ 
is $ (b!)^{j b^{k-1}}b^{- k}$.
\end{proof}

\section{Finite-size tests and perfect necklaces}

{\em ``Given a finite family of tests for randomness there is an infinite  sequence $x$ which passes all of them, 
but $x$ will be rejected by a new more refined test''}, proposed Norberto Fava to us. 
Our attempt to formalize this claim led to finite-size tests and perfect periodic sequences. 
The result is summarized in Proposition~\ref{kst}.

Let  $(X_0,X_1,\dots)$ be a sequence of random variables with values in a 
given alphabet~$\+A$ with at least two symbols. 
We say that the sequence is \emph{random} if the variables are uniformly distributed in~$\+A$ 
and mutually independent. To test if a sample $(x_0,\dots,x_{n-1})\in\+A^n$ comes from a 
random sequence we consider the following finite-size hypothesis testing~setup.
As usual, we write $\R$ for the set of real numbers.

\noindent(a) The \emph{hypothesis}
\begin{align*}
H_0&:\;\; (X_0,X_1,\dots) \hbox{ is random }
\end{align*}
\\
\noindent(b) A \emph{test-size} $k$ and a \emph{test function} $t:\+A^k\to \R$. Denote 
\[
\tau= E_0\big[t(X_0,\dots,X_{k-1})\big]=|\+A|^{-k}\sum_{(y_0,\dots,y_{k-1})\in\+A^k} t(y_0,\dots,y_{k-1}),
\] 
where $E_0$ is the expectation associated to the hypothesis $H_0$.
\\ \\
\noindent(c) A function $T_n:\+A^n\to\R$ defined by
\begin{align*}
T_n(x_0,\dots,x_{n-1}) = \Bigl|\frac1n\sum_{i=0}^{n-1} t(x_{i},\dots,x_{i+k-1})- \tau\Bigr|
\end{align*}
with periodic boundary conditions $x_{n+j}=x_j$. Thus, $T_n(x_0,\dots,x_{n-1})$ is the absolute 
difference between the empirical mean of $t$ for the sample and the expected value of $t$ under $H_0$.
\\ \\
\noindent(d) An \emph{error}  $\vep>0$ and the \emph{decision rule}
\begin{align*}
  \hbox{If  $T_n(x_0,\dots,x_{n-1})>\vep$ then reject the sample $(x_0,\dots,x_{n-1})$ as 
coming from $H_0$}. 
\end{align*}
In this case we say that \emph{the test $t$ rejects the sample $(x_0,\dots,x_{n-1})$}.

 This is called a test of size $k$ because rejection is decided as a function of the empirical mean of~$t$, 
a function of $k$ successive coordinates. Examples of finite-size tests include frequency test, 
block testing, number of runs in a block, longest run of ones in a block, etc. 
There are many (non-finite) tests, like the discrete Fourier transform test, 
the Kolmogorov-Smirnov test  and many others. 
Those tests also use some function $\tilde T_n$ of the sample, 
not necessarily based on the empirical mean of a $t$. 
The common feature is the use of the distribution of $\tilde T_n(X_1,\dots,X_n)$ 
under $H_0$ to compute the probability of rejection when $H_0$~holds. 

Tests for $H_0$ are used to check if a sequence of numbers 
produced by a random number generator can be considered random; see Knuth~\cite{MR3077153} and the  
battery of tests proposed by L'Ecuyer and Simard~\cite{MR2404400}. 
A nice account of the history of hypothesis testing is given by Lehmann~\cite{MR2798202}. 

In the usual hypotheses testing the sample-size $n$ is kept fixed. 
Assuming $H_0$ and repeating the test $j$ times with independent data, the proportion of times 
that the hypothesis is rejected converges as~$j\to\infty$ to the probability under~$H_0$ that 
$T_n(X_0,\dots,X_{n-1})>\vep$. Instead, we will take one infinite sequence, 
 test its first $n$ elements, record rejection for each~$n$ and take~$n\to\infty$. 

Let $x=(x_0,x_1, \dots)$ be an infinite sequence of symbols in~$\+A$. 
Fix a test-size $k$, a test-function $t$ of size $k$ and let $T_n$ be given by (c). 
We say that  \emph{$x$ passes the test $t$}~if 
\begin{align}
  \label{ht1} \tag{$*$}
\lim_{n\to\infty} T_n(x_0,\dots,x_{n-1}) = 0.
\end{align}
That is, for each $\vep>0$ there is an $n(x,\vep)$ such that for all $n>n(x,\vep)$ we~have  
\[
T_n(x_0,\dots,x_{n-1}) \le \vep. 
\]
In other words, fixing the test function $t$ of size $k$ and the error $\vep$, 
the test $t$ rejects $(x_0,\dots,x_{n-1})$ for at most a finite number of~$n$'s. 
When \eqref{ht1} 
does not hold we say that $t$ \emph{rejects~$x$}.

The random sequence $(X_0,X_1,\dots)$ of independently identically distributed uniform random variables in~$\+A$ 
passes any finite-size test $t$ almost surely. 
This is the same as saying that the set of real numbers in $[0,1]$ 
whose $|\+A|$-ary representation passes all finite tests has Lebesgue measure~$1$.

We say that the infinite sequence $x$ is \emph{$(k,m)$-perfect} if 
$x$ is periodic with period $m |\+A|^k$ and the necklace $[x_0\dots x_{m |\+A|^k-1}]$ is~$(k,m)$-perfect. 
Recall that $(k,1)$-perfect necklaces are exactly the de Bruijn necklaces of order $k$, 
so the following proposition considers infinite de Bruijn sequences of order $k$ as a special case:
if $x$ is de Bruijn of order $k$ there is a test of size $k+1$ that rejects~$x$.

\begin{Proposition}  \label{kst}
Assume alphabet  $\+A$ has at least two symbols. Let $m$ be a positive integer and 
let the infinite sequence  $x$ be $(k,m)$-perfect.  Then, the following holds:
\begin{enumerate}
\item  The infinite sequence $x$ passes every test of size $j\le k$. 

\item For each $h> k + \log_{|A|} m$ there exists a test $t$ of size $h$ 
                such that $t$ rejects~$x$.

\end{enumerate}
\end{Proposition}

\begin{proof} 
Let $b$ be the number of symbols in $\+A$. Thus, the period of $x$ has length $m b^k$.

1. Let $t$ be a test of size $k$.
For any positive integer $\ell$,  by periodiciy, 
\begin{align*}
T_{mb^k\ell} = \Bigl|\frac{1}{mb^k\ell} \sum_{i=0}^{mb^k\ell -1} t(x_i,\dots,x_{i+k-1})-\tau\Bigr| = 
\Bigl|\frac{\ell}{mb^k\ell} \sum_{i=0}^{mb^k-1} t(x_i,\dots,x_{i+k-1})-\tau\Bigr| = 0.
\end{align*}
because $x$ is $(k,m)$-perfect and the definition of $\tau$ in~(b).
Now take $j\in\{0,\dots, mb^k-1\}$ and use the above identity to get
  \[
(mb^k\ell+j)T_{mb^k\ell+j} = j\,T_{j} \le j \max|t-\tau| \le mb^k \max|t-\tau|,
\]
where $\max|t-\tau|=\max_{z_0,\dots,z_{k-1}}|t(z_0,\dots,z_{k-1})-\tau|$. 
Hence,
  \begin{align*}
T_{mb^k\ell+j} & \le \frac{mb^k}{mb^k\ell+j} \max|t-\tau|
\le \frac1\ell \max|t-\tau| \mathop{\longrightarrow}_{\ell\to\infty} 0.
\end{align*}
This shows that $x$ passes $t$.
Let $\tilde t$ be a test of size $j<k$. To see that $x$ also passes $\tilde t$  
define $t$ of size $k$~as
\[
 t(x_0,\dots,x_{k-1})= \tilde t(x_0,\dots,x_{j-1}).
\]

2.  Let $h$ be an integer such that  $h>k + \log_b m$.  
Then  $b^{h}> m b^k$ and there are more words  $w=w_0\dots w_{h-1}\in\+A^{h}$ than the possible  $m b^k$ 
places to start. Hence, there is at least one word $\tilde w$ of length $h$ not present in the sequence 
$x$ and  the test $t$ consisting on the indicator of $\tilde w$ rejects $x$. 
\end{proof}

\paragraph{Finite tests and normal numbers.}
As stated  by  Borel (see \cite{Bugeaud}), a real number is simply normal to base $b^k$ exactly when each 
block of length $k$ occurs in the $b$-ary expansion of $x$ with  asymptotic  frequency $b^{-k}$.
Hence, a real number is simply normal to base~$b^k$ if its $b$-ary expansion
 passes all tests up to size~$k$.
We have obtained that  for each $k$ and $b$, and for any~$m$,
 each $(k,m)$-perfect sequence in alphabet $\{0,1, \ldots, b-1\}$ 
is the $b$-ary expansion of a number that is simply normal  to base~$b^k$.
Borel defines   normality to base $b$ 
as simple normality to all bases~$b^k$, for every positive integer~$k$.
Henceforth, a number is normal to base~$b$ if its $b$-ary expansion passes all statistical tests of finite size.
Then, each instance of a number normal  to a given base provides an example of a sequence that  passes all finite-size tests.
Many  are known, such as \cite{Champernowne1933, Becherheiber2011} and the references in \cite{Bugeaud}.

\paragraph{Infinite tests and algorithmically random sequences.}
Martin L\"of introduced infinite tests defined in terms of computability~\cite{martinlof}.
These tests properly include all tests of finite size, so for every $k$ and $m$,
$(k,m)$-perfect sequences are rejected by these tests. 
The  infinite sequences that pass all these tests  are the  
 Martin L\"of random sequences, 
also known as the algorithmically random sequences.  
Due to the nature of the definition,  the algorithmically random sequences can not be computed
but some of them can be defined  at the first level of the Arithmetical Hierarchy~\cite{downey2010}.
\bigskip

\noindent
{\bf Acknowledgments.}
We thank Norberto Fava and Victor Yohai for motivating the question on the existence of 
periodic sequences that pass any finite family of  finite-size tests. 
We thank Liliana Forzani and Ricardo Fraiman for enlightening discussions. 

Alvarez and Becher are members of Laboratoire International Associ\'e INFINIS, 
Universit\'e Paris Diderot-CNRS / Universidad de Buenos Aires-CONICET).
Alvarez is supported by CONICET doctoral fellowship.
Becher  and Ferrari are  supported by the University of Buenos Aires and by CONICET.

\bibliographystyle{plain}

{\setstretch{0.9}

\bibliography{abfy}

\begin{thebibliography}{10}

\bibitem{Becherheiber2011}
Ver{\'o}nica Becher and Pablo~Ariel Heiber.
\newblock On extending de {B}ruijn sequences.
\newblock {\em Information Processing Letters}, 111(18):930--932, 2011.

\bibitem{berstel2007}
Jean Berstel and Dominique Perrin.
\newblock The origins of combinatorics on words.
\newblock {\em European Journal of Combinatorics}, 28(3):996--022, 2007.

\bibitem{Borel1909}
{{\'E}}mile Borel.
\newblock Les probabilit{\'e}s d{\'e}nombrables et leurs applications
  arithm{\'e}tiques.
\newblock {\em Supplemento Rendiconti del Circolo Matematico di Palermo},
  27:247--271, 1909.

\bibitem{Bugeaud}
Yann Bugeaud.
\newblock {\em Distribution modulo one and {D}iophantine approximation}, volume
  193 of {\em Cambridge Tracts in Mathematics}.
\newblock Cambridge University Press, Cambridge, 2012.

\bibitem{Champernowne1933}
David Champernowne.
\newblock The construction of decimals normal in the scale of ten.
\newblock {\em Journal of London Mathematical Society}, 3:254--260, 1933.

\bibitem{cvetovic}
Drago\v{s}~M. Cvetkovi\'{c}, Michael Doob, and Horst Sachs.
\newblock {\em Spectra of Graphs}.
\newblock Academic Press, 1980.

\bibitem{debruijn}
N.G. de~Bruijn.
\newblock A combinatorial problem.
\newblock {\em Indagationes Mathematicae}, 8:461--467, 1946.
\newblock Proc. Koninklijke Nederlandse Akademie v. Wetenschappen 49:758--764.

\bibitem{downey2010}
Rod Downey and Denis Hirschfeldt.
\newblock {\em Algorithmic Randomness and Complexity}.
\newblock Springer-Verlag New York, Inc., USA, 2010.

\bibitem{harary}
F.~Harary.
\newblock {\em Graph Theory}.
\newblock Addison-Wesley Publishing Co. Inc., Reading, Mass., 1969.

\bibitem{MR3077153}
Donald~E. Knuth.
\newblock {\em The art of computer programming. {V}ol. 2}.
\newblock Addison-Wesley, Reading, MA, 1998.
\newblock Seminumerical algorithms, Third edition [of MR0286318].

\bibitem{MR2404400}
Pierre L'Ecuyer and Richard Simard.
\newblock Test{U}01: a {C} library for empirical testing of random number
  generators.
\newblock {\em ACM Trans. Math. Software}, 33(4):Art. 22, 40, 2007.

\bibitem{MR2798202}
Erich~L. Lehmann.
\newblock {\em Fisher, {N}eyman, and the creation of classical statistics}.
\newblock Springer, New York, 2011.

\bibitem{lothaire1997}
M.~Lothaire.
\newblock {\em Combinatorics on words}.
\newblock Cambridge Mathematical Library. Cambridge University Press,
  Cambridge, 1997.

\bibitem{lothaire2002}
M.~Lothaire.
\newblock {\em Algebraic combinatorics on words}, volume~90 of {\em
  Encyclopedia of Mathematics and its Applications}.
\newblock Cambridge University Press, Cambridge, 2002.

\bibitem{martinlof}
Per Martin-L\"of.
\newblock The {D}efinition of {R}andom {S}equences.
\newblock {\em Information and Control}, 9(6):602--619, 1966.

\bibitem{tutte}
W.T. Tutte.
\newblock {\em Graph Theory}.
\newblock Addison-Wesley, 1984.

\end{thebibliography}

\begin{minipage}{\textwidth}
\noindent
Nicol\'as Alvarez
\\
Departamento de Ciencias e Ingenier\'ia de la Computaci\'on
\\
Universidad Nacional del Sur, Argentina.
\\
naa@cs.uns.edu.ar
\medskip\\
Ver\'onica Becher
\\
 Departmento de  Computaci\'on,   Facultad de Ciencias Exactas y Naturales
\\
Universidad de Buenos Aires \& CONICET, Argentina.
\\
vbecher@dc.uba.ar
\medskip\\
Pablo Ferrari
\\
Departmento de  Matem\'atica, Facultad de Ciencias Exactas y Naturales
\\
Universidad de Buenos Aires \& CONICET, Argentina.
\\
pferrari@dm.uba.ar
\medskip\\
Sergio Yuhjtman
\\
Departmento de  Matem\'atica, Facultad de Ciencias Exactas y Naturales\\
Universidad de Buenos Aires, Argentina.
\\
 syuhjtma@dm.uba.ar
\end{minipage}
}
\end{document}